\documentclass[12pt]{article}

\usepackage{latexsym, amssymb, amsmath, amscd, amsfonts, epsfig, graphicx, colordvi,verbatim,ifpdf}
\usepackage{amsfonts, amsmath, amssymb}
\usepackage{amssymb,amsfonts,amsmath,latexsym,epsfig,cite, psfrag,eepic,color}
\usepackage{amscd,graphics}
\usepackage{latexsym, amssymb,  amsmath,amscd, amsfonts, epsfig, graphicx, colordvi,amsthm}

\usepackage{graphicx}
\usepackage{color}
\usepackage{ifpdf}
\usepackage{fancybox}

\newtheorem{thm}{Theorem}[section]

\newtheorem{defi}[thm]{Definition}

\def\pf{\noindent{\it Proof.} }
\def\qed{\nopagebreak\hfill{\rule{4pt}{7pt}}
\medbreak}

\numberwithin{equation}{section}

\def\qed{\nopagebreak\hfill{\rule{4pt}{7pt}}
\medbreak}

\newlength{\boxedparwidth}
  {\begin{center} \begin{tabular}{|@{\hspace{.315in}}c@{\hspace{.15in}}|}
                  \hline \\ \begin{minipage}[t]{\boxedparwidth}
                  \setlength{\parindent}{.25in}}%
  {\end{minipage} \\ \\ \hline \end{tabular} \end{center}}

\begin{document}
\begin{center}

{ \large\bf  On the Number of Partitions with Designated Summands}
\end{center}

\vskip 5mm

\begin{center}
{  William Y.C. Chen}$^{1}$,    {  Kathy Q. Ji}$^{2}$, { Hai-Tao
Jin}$^{3}$ and {   Erin  Y.Y. Shen}$^{4}$ \vskip 2mm

   Center for Combinatorics, LPMC-TJKLC\\
   Nankai University, Tianjin 300071, P.R. China

   \vskip 2mm

    $^1$chen@nankai.edu.cn, $^2$ji@nankai.edu.cn, $^3$jinht1006@mail.nankai.edu.cn, $^4$shenyiying@mail.nankai.edu.cn
\end{center}

\vskip 6mm \noindent {\bf Abstract.} Andrews, Lewis and Lovejoy introduced the
partition function $PD(n)$ as the number of partitions of $n$ with designated summands, where we assume that
among parts with equal size, exactly one  is designated.
They proved that   $PD(3n+2)$ is divisible by $3$.
We obtain a Ramanujan type identity for the
generating function of $PD(3n+2)$ which implies the congruence
of Andrews, Lewis and Lovejoy. For $PD(3n)$, Andrews, Lewis and Lovejoy showed that the  generating function   can be expressed as an infinite product of powers of $(1-q^{2n+1})$ times  a function $F(q^2)$.  We find an explicit formula for $F(q^2)$,
which leads to a formula for the generating function of $PD(3n)$.
We also obtain a formula for the generating function of $PD(3n+1)$.
 Our proofs rely on Chan's identity on Ramanujan's cubic continued fraction and some identities on cubic theta functions.  By introducing   a rank for the partitions with designed summands,
 we give a combinatorial interpretation of the
 congruence  of Andrews, Lewis and Lovejoy.

\section{Introduction}

 Andrews, Lewis and Lovejoy \cite{Andrews-Lewis-Lovejoy-2002}
 investigated the number of
 partitions with designated summands which are defined on ordinary partitions by designating exactly one part among parts with equal size.
 Let $PD(n)$ denote the number of partitions of $n$ with designated summands.
  For example, there are  ten partitions of $4$ with designated summands:\emph{}
\begin{align*}
\begin{array}{ccccc}
4',      & 3'+1',  & 2'+2,    & 2+2',    & 2'+1'+1,  \\[5pt]
2'+1+1', &1'+1+1+1,& 1+1'+1+1,& 1+1+1'+1,& 1+1+1+1'.
\end{array}
\end{align*}

The notion of partitions with designated summands goes back to
 MacMahon \cite{MacMahon-1919}. He considered
 partitions with designated summands and with
  exactly $k$ different sizes, see also Andrews and Rose  \cite{Andrews-Rose}.
   Andrews, Lewis and Lovejoy \cite{Andrews-Lewis-Lovejoy-2002}  derived  the following generating function of $PD(n)$.

\begin{thm}  We have
\begin{align}\label{PD(n)}
\sum_{n=0}^{\infty}PD(n)q^n&=\frac{(q^6;q^6)_{\infty}}{(q;q)_{\infty}
(q^2;q^2)_{\infty}(q^3;q^3)_{\infty}},
\end{align}
where $(a;q)_\infty$ stands for the $q$-shifted factorial
\[ (a;q)_{\infty}=\prod_{n=1}^{\infty}(1-aq^{n-1}),\quad     |q|<1. \]
\end{thm}

By using modular forms and $q$-series identities,  Andrews, Lewis and Lovejoy showed that the partition function $PD(n)$ has many interesting divisibility properties. In particular, they obtained the following Ramanujan type congruence.

\begin{thm}{\rm (\!\!\cite[Corollary 7]{Andrews-Lewis-Lovejoy-2002})}\label{all-1} For $n\geq 0$, we have
\begin{align}
\label{modPD(3n+2)}PD(3n+2)\equiv 0\pmod{3}.
\end{align}
\end{thm}

In this paper, we obtain the following Ramanujan type identity for the
generating function of $PD(3n+2)$ which implies the above congruence.

\begin{thm}\label{PD(3n+2)-tm}We have
\begin{equation}\label{PD(3n+2)}
\sum_{n=0}^{\infty}PD(3n+2)q^n=3\frac{(q^3;q^6)_{\infty}^{3}
(q^6;q^6)_{\infty}^{6}}{(q;q^2)_{\infty}^{5}(q^2;q^2)_{\infty}^{8}}.
\end{equation}
\end{thm}

 Andrews, Lewis and Lovejoy also obtained  explicit formulas for    the generating functions for $PD(2n)$ and $PD(2n+1)$ by using Euler's algorithm for infinite products \cite[P. 98]{Andrews-1976} and Sturm's criterion \cite{Sturm-1984}. As for $PD(3n)$, they showed that
    the generating function  permits the following  form.

\begin{thm}{\rm (\!\!\cite[Theorem 23]{Andrews-Lewis-Lovejoy-2002})}\label{all-2}
Define $c(n)$ uniquely by
\begin{align}
\sum_{n=0}^{\infty}PD(3n)q^n=\prod_{n=1}^{\infty}(1-q^n)^{-c(n)},
\end{align}
then for all positive $n$,
\begin{align*}
c(6n+1)=5,\\
c(6n+3)=2,\\
c(6n+5)=5.
\end{align*}
\end{thm}

Equivalently, the above theorem says that  there exists a  series $F(q^2)$ such that

\begin{align}\label{Fq2}
\sum_{n=0}^{\infty}PD(3n)q^n&=\frac{1}{(q;q^6)_{\infty}^{5}(q^3;q^6)_{\infty}^{2}
(q^5;q^6)_{\infty}^{5}}\times F(q^2).
\end{align}

 In this paper, we find an explicit formula for $F(q^2)$, that is,
 \begin{align}
F(q^2)=\frac{(q^4;q^4)_{\infty}^{6}(q^6;q^6)_{\infty}^{4}}{(q^2;q^2)_{\infty}^{10}(q^{12};q^{12})_{\infty}^{2}}
                                +3q^2\frac{(q^{12};q^{12})_{\infty}^{6}}{(q^2;q^2)_{\infty}^{6}(q^4;q^4)_{\infty}^{2}},
\end{align}
which leads to the following generating function  of $PD(3n)$.

\begin{thm}\label{PD(3n)-tm}We have
\begin{align}
\label{PD(3n)}&\sum_{n=0}^{\infty}PD(3n)q^n
=\frac{1}{(q;q^6)_{\infty}^{5}(q^3;q^6)_{\infty}^{2}(q^5;q^6)_{\infty}^{5}}\times \nonumber \\
                               & \hskip 4cm \left(\frac{(q^4;q^4)_{\infty}^{6}(q^6;q^6)_{\infty}^{4}}{(q^2;q^2)_{\infty}^{10}(q^{12};q^{12})_{\infty}^{2}}
                                +3q^2\frac{(q^{12};q^{12})_{\infty}^{6}}{(q^2;q^2)_{\infty}^{6}(q^4;q^4)_{\infty}^{2}} \right).
\end{align}
\end{thm}

In fact, we obtain explicit formulas for the $3$-dissection of the generating
function    of $PD(n)$, which include the following generating function
for $PD(3n+1)$.

\begin{thm}\label{PD(3n+1)-tm}We have
\begin{equation}\label{PD(3n+1)}
\sum_{n=0}^{\infty}PD(3n+1)q^n
=\frac{(q^3;q^6)_{\infty}^{3}(q^6;q^6)_{\infty}^{6}}
{(q;q^2)_{\infty}^{5}(q^2;q^2)_{\infty}^{8}}
                                \left(4q\frac{(q;q^2)_{\infty}^{2}}
                                {(q^3;q^6)_{\infty}^{6}}+\frac{(q^3;q^6)_{\infty}^{3}}
                                {(q;q^2)_{\infty}}\right).
\end{equation}
\end{thm}

Our dissection formulas rely on the Chan's identity  on Ramanujan's cubic continued fraction \cite{Chan-2008}  and cubic theta functions \cite{Hirschhorn-Garvan-Borwein-1993, Berndt-Bhargava-Garvan-1995}.  In Section 3, we shall give a combinatorial interpretation of  the congruence $PD(3n+2) \equiv 0 \pmod 3$ by introducing a rank for the partitions with designed summands.

\section{Proofs}

In this section, we give proofs of  the generating functions for $PD(3n),$ $PD(3n+1)$ and $PD(3n+2)$  by employing Chan's
identity on Ramanujan's cubic continued fraction.
 It should be noted that the generating function
 of $PD(3n)$ derived this way does not directly imply a formula
 for $F(q^2)$. To compute $F(q^2)$, we shall make use of
 some identities on cubic theta functions.

Recall that   Ramanujan's cubic continued fraction $\upsilon(q)$ is given by
\[\upsilon(q):=\frac{q^{\frac{1}{3}}}{1}_{\ +\ }\frac{q+q^2}{1}
_{\ +\ }\frac{q^2+q^4}{1}_{\ +\ \ldots}.
\]
It is  known that
\[\upsilon(q)=q^{\frac{1}{3}}\frac{(q;q^2)_{\infty}}{(q^3;q^6)_{\infty}^3},\]
see Andrews and Berndt \cite[P. 94]{Andrews-Berndt-2005}.
 The following  identity is due to Chan and  will be used in our derivation of the 3-dissection formulas.

\begin{thm}{\rm (\!\!\cite[Eq. (13)]{Chan-2008})}\label{chan}
We have
\begin{align}\label{chan-id}
\frac{1}{(q;q)_{\infty}(q^2;q^2)_{\infty}}
&=\frac{(q^9;q^9)_{\infty}^3(q^{18};q^{18})
_{\infty}^3}{(q^3;q^3)_{\infty}^4
(q^6;q^6)_{\infty}^4} \times \nonumber \\[2pt]
&\qquad \left\{\left(\frac{1}{x^2(q^3)}-2q^3x(q^3)\right)
+q\left(\frac{1}{x(q^3)}+4q^3x^2(q^3)\right)+3q^2\right\},
\end{align}
where
\[x(q)=q^{-\frac{1}{3}}v(q)=\frac{(q;q^2)_{\infty}}{(q^3;q^6)_{\infty}^3}.\]
 \end{thm}

\noindent{\it Proof of Theorems \ref{PD(3n+2)-tm} and \ref{PD(3n+1)-tm}. }
Multiplying  both sides of \eqref{chan-id} by
\[\frac{(q^6;q^6)_{\infty}}{(q^3;q^3)_{\infty}},\]
 we find
\begin{align}\label{PD-3dis}
\frac{(q^6;q^6)_\infty}{(q;q)_{\infty}(q^2;q^2)_{\infty}(q^3;q^3)_\infty}
&=\frac{(q^9;q^9)_{\infty}^3(q^{18};q^{18})
_{\infty}^3}{(q^3;q^3)_{\infty}^5
(q^6;q^6)_{\infty}^3}\left\{\left(\frac{1}{x^2(q^3)}-2q^3x(q^3)\right)\right.
\nonumber \\[2pt]
&\qquad \qquad \qquad \qquad \qquad \left.+q\left(\frac{1}{x(q^3)}+4q^3x^2(q^3)\right)+3q^2\right\}.
\end{align}
Observe that the left-hand side of \eqref{PD-3dis} is the generating function for $PD(n)$. Extracting those terms  involving the powers $q^{3n},\,q^{3n+1}$ and $q^{3n+2}$, respectively, we deduce that
\begin{align}
\sum_{n=0}^{\infty}PD(3n)q^{3n}
\quad &=\frac{(q^9;q^9)_{\infty}^3(q^{18};q^{18})_{\infty}^3}
{(q^3;q^3)_{\infty}^5(q^6;q^6)_{\infty}^3}
\left(-2q^3x(q^3)+\frac{1}{x^2(q^3)}\right),\label{PD(3n)temp}\\
\sum_{n=0}^{\infty}PD(3n+1)q^{3n+1}
&=q\frac{(q^9;q^9)_{\infty}^3(q^{18};q^{18})_{\infty}^3}
{(q^3;q^3)_{\infty}^5(q^6;q^6)_{\infty}^3}                                              \left(4q^3x^2(q^3)+\frac{1}{x(q^3)}\right), \label{PD(3n+1)temp}\\
\sum_{n=0}^{\infty}PD(3n+2)q^{3n+2}
&=3q^2\frac{(q^9;q^9)_{\infty}^3(q^{18};q^{18})_{\infty}^3}
{(q^3;q^3)_{\infty}^5(q^6;q^6)_{\infty}^3}.\label{PD(3n+2)temp}
\end{align}
Thus Theorem \ref{PD(3n+2)-tm} can be deduced from \eqref{PD(3n+2)temp} by dividing both sides   by $q^2$ and substituting $q^3$ by $q$.
Similarly, Theorem \ref{PD(3n+1)-tm} can be deduced from \eqref{PD(3n+1)temp} by dividing both sides  by $q$ and substituting $q^3$ by $q$. This completes the proof. \qed

 It turns out that $F(q^2)$ can be computed  with the aid of some identities for cubic theta functions.    These functions are introduced by Borwein,  Borwein and Garvan \cite{Borwein1994} and are defined by \begin{align*}
a(q)&=\sum_{m,n=-\infty}^{\infty} q^{m^2+mn+n^2},\\
b(q)&=\sum_{m,n=-\infty}^{\infty} \omega^{m-n} q^{m^2+mn+n^2}, \quad \omega=e^{2\pi i/3},\\
c(q)&=\sum_{m,n=-\infty}^{\infty} q^{m^2+mn+n^2+m+n}.
\end{align*}
Recall that
 \begin{equation}
 c(q)=3\frac{(q^3;q^3)_{\infty}^3}{(q;q)_{\infty}},
\end{equation}
see Berndt, Bhargava and Garvan \cite[Eq. (5.5)]{Berndt-Bhargava-Garvan-1995}.
We shall also use the following identities for $a(q)$ and $c(q)$
\begin{eqnarray}
a(q)&=&a(q^4)+6q\frac{(q^4;q^4)_{\infty}^{2}(q^{12};q^{12})_{\infty}^{2}}
{(q^2;q^2)_{\infty}(q^6;q^6)_{\infty}},\label{aa}\\
c(q)&=&qc(q^4)+3\frac{(q^4;q^4)_{\infty}^{3}(q^6;q^6)_{\infty}^{2}}
{(q^2;q^2)_{\infty}^{2}
(q^{12};q^{12})_{\infty}},\label{cc}\\
a(q)&=&a(q^2)+2q\frac{c^2(q^2)}{c(q)}.\label{ca}
\end{eqnarray}
Identity \eqref{aa} for $a(q)$  and identity  \eqref{cc} for $c(q)$ are due to Hirschhorn, Garvan, and Borwein \cite[Eqs.(1.36) and (1.34)]{Hirschhorn-Garvan-Borwein-1993}.
Identity  \eqref{ca} for $a(q)$ and $c(q)$  is obtained by Berndt, Bhargava, Garvan \cite[Eq. (6.3)]{Berndt-Bhargava-Garvan-1995}.

We obtain the following  identity on Ramanujan's cubic continued fraction.

\begin{thm}\label{dissect-2}Let
\[x(q)=q^{-\frac{1}{3}}v(q)=\frac{(q;q^2)_{\infty}}{(q^3;q^6)_{\infty}^3}.\]
We have
\begin{align}\label{dissect-eq}
\frac{1}{x^2(q)}-2qx(q)=3q^2\frac{(q^2;q^2)_{\infty}^{2}(q^{12};q^{12})_{\infty}^{6}}{(q^4;q^4)_{\infty}^{2}(q^6;q^6)_{\infty}^{6}}
       +\frac{(q^4;q^4)_{\infty}^{6}}{(q^2;q^2)_{\infty}^{2}(q^6;q^6)_{\infty}^{2}(q^{12};q^{12})_{\infty}^{2}}.
\end{align}
\end{thm}

\pf We first establish a connection between Ramanujan's cubic continued fraction $\upsilon(q)$  and the cubic theta function $c(q)$. It is easy to check that
\begin{align}
\label{A} \frac{1}{x^2(q)}&=\frac{(q^3;q^6)_{\infty}^{6}}{(q;q^2)_{\infty}^{2}}
            =\frac{(q^2;q^2)_{\infty}^{2}}{(q^6;q^6)_{\infty}^{6}}\times\left(\frac{(q^3;q^3)_{\infty}^{3}}{(q;q)_{\infty}}\right)^2
            =\frac{(q^2;q^2)_{\infty}^{2}}{9(q^6;q^6)_{\infty}^{6}}\times c^2(q),\\
\label{B} 2qx(q)&=2q\frac{(q;q^2)_{\infty}}{(q^3;q^6)_{\infty}^{3}}
            =2q\frac{(q^6;q^6)_\infty^3}{(q^2;q^2)_\infty}\times
            \left(\frac{(q;q)_\infty}{(q^3;q^3)_\infty^3}\right)
            =6q\frac{(q^6;q^6)_{\infty}^{3}}{(q^2;q^2)_{\infty}}\times \frac{1}{c(q)}.
\end{align}

We now consider  the 2-dissection of $ 1/x^2(q)$.
 Identity \eqref{cc}  can be viewed as the $2$-dissection of $c(q)$. Hence
 we deduce that
  \begin{align*}
c^2(q)=\left(q^2c^2(q^4)+9\frac{(q^4;q^4)_{\infty}^{6}(q^6;q^6)_{\infty}^{4}}
{(q^2;q^2)_{\infty}^{4}(q^{12};q^{12})_{\infty}^{2}}\right)
+q\left(6c(q^4)\frac{(q^4;q^4)_{\infty}^{3}(q^6;q^6)_{\infty}^{2}}
{(q^2;q^2)_{\infty}^{2}(q^{12};q^{12})_{\infty}}\right).
\end{align*}
This yields the 2-dissection of $1/x^2(q)$,
\begin{align}\label{A123}
  \frac{1}{x^2(q) } &=\left(q^2c^2(q^4)\frac{(q^2;q^2)_{\infty}^{2}}{9(q^6;q^6)_{\infty}^{6}}+\frac{(q^2;q^2)_{\infty}^{2}}{(q^6;q^6)_{\infty}^{6}}
    \frac{(q^4;q^4)_{\infty}^{6}(q^6;q^6)_{\infty}^{4}}
    {(q^2;q^2)_{\infty}^{4}(q^{12};q^{12})_{\infty}^{2}} \right)\nonumber\\[3pt]
    &\qquad
    +q\left(6c(q^4)\frac{(q^2;q^2)_{\infty}^{2}}{9(q^6;q^6)_{\infty}^{6}} \frac{(q^4;q^4)_{\infty}^{3}(q^6;q^6)_{\infty}^{2}}
    {(q^2;q^2)_{\infty}^{2}(q^{12};q^{12})_{\infty}}\right)\nonumber\\[7pt]
    &= \left( q^2\frac{(q^2;q^2)_{\infty}^{2}(q^{12};q^{12})_{\infty}^{6}}
    {(q^4;q^4)_{\infty}^{2}(q^6;q^6)_{\infty}^{6}}
    +\frac{(q^4;q^4)_{\infty}^{6}}{(q^2;q^2)_{\infty}^{2}
    (q^6;q^6)_{\infty}^{2}(q^{12};q^{12})_{\infty}^{2}}\right)\nonumber\\[3pt]
    &\qquad + 2q\left(\frac{(q^4;q^4)_{\infty}^{2}(q^{12};q^{12})_{\infty}^{2}}
    {(q^6;q^6)_{\infty}^{4}}\right).
\end{align}

 Next, we aim to  derive the 2-dissection of $q/c(q)$. By \eqref{ca}, we find
\begin{equation}\label{temp}
\frac{q}{c (q)}=\frac{a(q)-a(q^2)}{2c^2(q^2)}.
\end{equation}
Substituting  \eqref{aa} into \eqref{temp}, we arrive at
\begin{equation}\label{temp-2}
\frac{q}{c (q)}=\frac{1}{2c^2(q^2)}\left(a(q^4)+6q\frac{(q^4;q^4)_{\infty}^{2}
(q^{12};q^{12})_{\infty}^{2}}{(q^2;q^2)_{\infty}
(q^6;q^6)_{\infty}}-a(q^2)\right).
 \end{equation}
Using \eqref{ca} with $q$ replaced by $q^2$, we get
\[  a(q^2)-a(q^4)= 2q^2\frac{c^2(q^4)}{c(q^2)} .\]
Hence  \eqref{temp-2} can be written as
  \begin{align*}
  \frac{q}{c (q)}&=\frac{1}{2c^2(q^2)}\left(-2q^2\frac{c^2(q^4)}{c(q^2)}
  +6q\frac{(q^4;q^4)_{\infty}^{2}(q^{12};q^{12})_{\infty}^{2}}
  {(q^2;q^2)_{\infty}(q^6;q^6)_{\infty}}\right)\\[6pt]
            &=-q^2\frac{c^2(q^4)}{c^3(q^2)}+3q
\frac{(q^4;q^4)_{\infty}^{2}(q^{12};q^{12})_{\infty}^{2}}{c^2(q^2)(q^2;q^2)_{\infty}(q^6;q^6)_{\infty}}.
\end{align*}
Thus, we obtain  the following 2-dissection of $2qx(q)$,
\begin{align}\label{2qB12}
2qx(q)&=-6q^2\frac{ (q^6;q^6)_{\infty}^{3}c^2(q^4)}{(q^2;q^2)_{\infty}c^3(q^2)}
      +18q
      \frac{(q^6;q^6)_{\infty}^{3}(q^4;q^4)_{\infty}^{2}(q^{12};q^{12})_{\infty}
      ^{2}}{c^2(q^2)(q^2;q^2)_{\infty}^2(q^6;q^6)_{\infty}}\nonumber \\[6pt]
    &=-2q^2\frac{(q^2;q^2)_{\infty}^{2}(q^{12};q^{12})_{\infty}^{6}}{(q^6;q^6)_{\infty}^{6}(q^4;q^4)_{\infty}^{2}}
     +2q\frac{(q^4;q^4)_{\infty}^{2}(q^{12};q^{12})_{\infty}^{2}}
     {(q^6;q^6)_{\infty}^{4}}.
\end{align}
Subtracting \eqref{2qB12}  from \eqref{A123}, we obtain \eqref{dissect-eq}.  This completes the proof.
 \qed

\noindent {\it Proof of Theorem \ref{PD(3n)-tm}.} Substituting $q^3$ with  $q$
in \eqref{PD(3n)temp}, we obtain
\begin{align}
\sum_{n=0}^{\infty}PD(3n)q^{n}
 &=\frac{(q^3;q^3)_{\infty}^3(q^{6};q^{6})_{\infty}^3}
{(q;q)_{\infty}^5(q^2;q^2)_{\infty}^3}
\left(-2qx(q)+\frac{1}{x^2(q)}\right)\nonumber \\[4pt]
&=\frac{ (q^3;q^6)_\infty^3(q^6;q^6)_{\infty}^6}{(q;q^2)_{\infty}^5 (q^2;q^2)_{\infty}^8}\left(-2qx(q)+\frac{1}{x^2(q)}\right) \nonumber \\[4pt]
&=\frac{1}{(q;q^6)_{\infty}^{5}(q^3;q^6)_{\infty}^{2}(q^5;q^6)_{\infty}^{5}}\times\frac{  (q^6;q^6)_{\infty}^6}{ (q^2;q^2)_{\infty}^8}\left(-2qx(q)+\frac{1}{x^2(q)}\right). \label{temp-3}
\end{align}
Applying \eqref{dissect-eq} to \eqref{temp-3}, we are led to the generating function for $PD(3n)$ in Theorem \ref{PD(3n)-tm}. This completes the proof. \qed

\section{A combinatorial interpretation}

In this section, we  give a combinatorial interpretation of the
 congruence $PD(3n+2) \equiv 0 \pmod 3$. In doing so,
  we introduce a rank for partitions with designated summands.
 We call this rank the $pd$-rank which enables us to divide the
   set of   partitions of $3n+2$ with designated summands   into three equinumerous classes.  The definition of the $pd$-rank is based on the
   following representation of a partition with designated summands by
   a pair of partitions.

\begin{thm}\label{lambdaab}  There is a
bijection $\Delta$ between the set of partitions of $n$ with designed summands  and the set of pairs  of partitions  $(\alpha, \beta)$   of $n$,
where $\alpha$ is an ordinary partition and $\beta$  is a partition into parts $\not\equiv
\pm1 \pmod{6}$.
\end{thm}

To give a  proof of the above theorem, we shall use
the bijective proof of the following theorem of MacMahon given by
  Andrews, Eriksson, Petrov and Romik \cite{Andrews-Eriksson-Petrov-Romik-2007}.

\begin{thm}{\rm (\!\!\cite{MacMahon-1960})}\label{MacMahonthm}
The number of partitions of an integer $n$ into parts $\not\equiv
\pm1 \pmod{6}$ equals the number of partitions of $n$ not containing any part exactly once.
\end{thm}

\noindent{\it Proof of Andrews, Eriksson, Petrov and Romik.}
We construct a bijection $\Phi$ from the set $\mathcal{C}_n$ of partitions of $n$ not containing any part exactly once to the set $\mathcal{B}_n$ of partitions of $n$ into parts not congruent to $\pm1$ mod $6$. To describe
 the map $\Phi$, let $\mu$ be a partition in $\mathcal{C}_n$.
 Write  $\mu$ as in the form of $(1^{m_1}2^{m_2}\cdots l^{m_l})$, where $m_k$ is the multiplicity of $k$ so that $n=\sum_{k=1}^{t}km_k$. Since $m_k \not= 1$ for any $k$,
   there is a unique way to write $m_k$  as $m_k=s_k+t_k$, where $s_k\in \{0,3\}$ and $t_k\in\{0,2,4,6,8,\ldots\}$. Now, the partition $\lambda=\Phi(\mu)=(1^{b_1}2^{b_2}\ldots)$ is determined as follows:
\begin{align*}
&b_{6k+1}=0, \hskip 3cm  b_{6k+5}=0,\\[2pt]
&b_{6k+2}=\frac{1}{2}t_{3k+1},\hskip 2.2cm  b_{6k+4}=\frac{1}{2}t_{3k+2},\\
&b_{6k+3}=\frac{1}{3}s_{2k+1}+t_{6k+3},\hskip 0.8cm  b_{6k+6}=\frac{1}{3}s_{2k+2}+t_{6k+6}.
\end{align*}
It is evident that $\lambda$ is a partition into parts not congruent to $\pm1$ mod $6$. It is also apparent that one can recover the partition $\mu$ from
$\lambda$ by reversing the above procedure. Hence $\Phi$ is a bijection.
This completes the proof. \qed

 We are now in a position to present the  proof of Theorem \ref{lambdaab} by using the bijection $\Phi$. \medskip

 \noindent {\it Proof of Theorem \ref{lambdaab}.} Let $\lambda$ be a partition of $n$ with designed summands. We wish to construct a pair  of partitions $(\alpha,\beta)$ of $n$, where $\alpha$ is an ordinary partition and $\beta$ is a partition  into parts $\not\equiv \pm1 \pmod{6}$.

 Suppose $t$ is a magnitude that appears in $\lambda$ and there are $m_t$ parts equal to $t$  among which  the $i$-th part is designated. There are two cases.
\begin{itemize}
\item
  If $i=1$, then move all the parts equal to $t$ (including the designated part) in $\lambda$ to the partition $\alpha$.

\item
If $i\neq 1$, then move  $i$  parts equal to $t$ in $\lambda$ to $\gamma$  and $(m_t-i)$  parts equal to $t$ in $\lambda$  to $\alpha$.
\end{itemize}
It can be seen that each part in $\gamma$   occurs at least twice. Let $\beta=\Phi(\gamma)$. It is clear    that  $\beta$ is a partition  into parts $\not\equiv
\pm1 \pmod{6}$ and the above procedure can be reversed. Hence $\Delta$ is a bijection. This completes the proof.  \qed

  The $pd$-rank of a partition $\lambda$ with designated summands can
  be defined in terms of the pair of partitions  $(\alpha, \beta)$ under the
  map  $\Delta$.

\begin{defi}
Let $\lambda$ be a partition with designated summands and let $(\alpha,\beta)=\Delta(\lambda)$.   Then the $pd$-rank of  $\lambda$,  denoted  $r_d(\lambda)$, is defined by
\begin{align}
r_d(\lambda)=l_e(\alpha)-l_e(\beta),
\end{align}
where
$l_e(\alpha)$ is the number of even parts of $\alpha$ and
and $l_e(\beta)$ is the number of even parts of $\beta$.
\end{defi}

The following theorem  shows that the  $pd$-rank  can be
  used to divide the set of   partitions of $3n+2$ with designated summands   into three equinumerous classes.

\begin{thm}\label{thm-c}
For $i=0,1,2$, let $N_d(i,3;n)$ denote the number of partitions of $n$ with designated
summands with $pd$-rank congruent to $i \pmod{3}$. Then we have
\begin{align}
N_d(0,3;3n+2)=N_d(1,3;3n+2)=N_d(2,3;3n+2).
\end{align}
\end{thm}

\pf  Let  $N_d(m;n)$ denote the number of partitions of $n$ with designated
summands with $pd$-rank $m$. By the definition of the $pd$-rank, we see that
\begin{align}
\sum_{n=0}^{\infty}\sum_{m=-\infty}^{\infty}N_d(m;n)z^mq^n=\frac{ 1}{(zq^2;q^2)_{\infty}(q;q^2)_\infty}\times \frac{1}{(z^{-1}q^2;q^2)_{\infty}(q^3;q^6)_{\infty}}.
\end{align}
Setting $z=\zeta=e^{\frac{2\pi i}{3}}$, we find that
\begin{align}
\sum_{n=0}^{\infty}\sum_{m=-\infty}^{\infty}N_d(m;n)\zeta^mq^n
&=\sum_{n=0}^{\infty}\sum_{i=0}^2N_d(i,3;n)\zeta^i q^n \nonumber  \\
&=\frac{1}{(\zeta
    q^2;q^2)_{\infty}(q;q^2)_\infty(\zeta^{-1}q^2;q^2)_{\infty}(q^3;q^6)_\infty} \nonumber \\
  &=\frac{(-q^3;q^3)_{\infty}}{(q;q^2)_\infty(\zeta
    q^2;q^2)_{\infty}(\zeta^{-1}q^2;q^2)_{\infty}}. \label{temp-c}
\end{align}
Multiplying the right hand side of \eqref{temp-c}  by
\[\frac{(q^2;q^2)_{\infty}}{(q^2;q^2)_{\infty}},\]
and noting that
\[(1-x)(1-x\zeta)(1-x\zeta^2)=1-x^3,
\]
we deduce that
\begin{align*}
 \sum_{n=0}^{\infty}\sum_{i=0}^2N_d(i,3;n)\zeta^i q^n&=\frac{(-q^3;q^3)_{\infty}}{(q;q^2)_{\infty}(\zeta q^2;q^2)_{\infty}(\zeta^{-1}q^2;q^2)_{\infty}}\times
     \frac{(q^2;q^2)_{\infty}}{(q^2;q^2)_{\infty}}\\
     &=\frac{(q^2;q^2)_{\infty}}{(q;q^2)_{\infty}}\times
     \frac{(-q^3;q^3)_{\infty}}{(q^6;q^6)_{\infty}}.
\end{align*}
By Gauss's identity \cite[P. 23]{Andrews-1976}
\[\frac{(q^2;q^2)_{\infty}}{(q;q^2)_{\infty}}=\sum_{n=0}^{\infty}q^{{n+1 \choose
      2}},
\]
we get
\begin{align}\label{temp-c-2}
\sum_{n=0}^{\infty}\sum_{i=0}^2N_d(i,3;n)\zeta^i q^n=\frac{(-q^3;q^3)_{\infty}}{(q^6;q^6)_{\infty}}\sum_{n=0}^{\infty}q^{{n+1 \choose
      2}}.
\end{align}
Since
\[{n+1 \choose
      2}\equiv 0 \  \text{ or } \    1 \pmod 3,\]
       the coefficient of
$q^{3n+2}$ in   \eqref{temp-c-2} is zero.
It follows that
\[N_d(0,3;3n+2)+N_d(1,3;3n+2)\zeta+N_d(1,3;3n+2)\zeta^2=0. \]
Since $1+\zeta+\zeta^2$ is the minimal polynomial in $\mathbb{Z}[\zeta]$,  we conclude that
\[
N_d(0,3;3n+2)=N_d(1,3;3n+2)=N_d(2,3;3n+2).
\] This completes the proof.\qed

For example, for $n=5$, we have $PD(5)=15$. The fifteen partitions
of $5$ with designated summands, the corresponding pairs of
partitions, along with the $pd$-ranks modulo 3 are listed in  Table 3.1.
It can be checked that
\[ N_d(0,3;5)=N_d(1,3;5)=N_d(2,3;5)=5.\]

\begin{table}[h]\centering
\begin{tabular}  {c|c|c}
$\lambda$               & $(\alpha,\beta)=\Delta(\lambda)$      & $r_d(\lambda) \pmod 3$\\[4pt]
\hline
$5'$                    &  $(5,\ \emptyset)$           & $0$\\[4pt]
\hline
$4'+1'$             &  $(4+1,\ \emptyset)$       & $1$\\[4pt]
\hline
$3'+2'$            &  $(3+2,\  \emptyset)$       & $1$\\[4pt]
\hline
$3'+1'+1$           &  $(3+1+1,\ \emptyset)$     & $0$\\[4pt]
\hline
 $3'+1+1'$                &  $(3,\  2)$               & $2$\\[4pt]
\hline
$2'+2+1'$         &  $(2+2+1,\ \emptyset)$     & $2$\\[4pt]
\hline
$2+2'+1'$              &  $(1,\ 4)$               & $2$\\[4pt]
\hline
$2'+1'+1+1$     &  $(2+1+1+1,\ \emptyset)$   & $1$\\[4pt]
\hline
$2'+1+1'+1$              &  $(2+1,\ 2)$             & $0$\\[4pt]
\hline
$2'+1+1+1'$           &  $(2,\ 3)$               & $1$\\[4pt]
\hline
$1'+1+1+1+1$  &  $(1+1+1+1+1,\ \emptyset)$  & $0$\\[4pt]
\hline
$1+1'+1+1+1$        &  $(1+1+1,\ 2)$           & $2$\\[4pt]
\hline
$1+1+1'+1+1$          &  $(1+1,\  3)$             & $0$\\[4pt]
\hline
$1+1+1+1'+1$        &  $(1,\ 2+2)$             & $1$\\[4pt]
\hline
$1+1+1+1+1'$   &  $(\emptyset,\ 3+2)$       & $2$
\end{tabular}
 \caption{ The case for $n=5$.}
\end{table}

 \vspace{.2cm}

\noindent{\bf Acknowledgments.} This work was supported by the 973
Project, the PCSIRT Project of the Ministry of Education, and the
National Science Foundation of China.

\end{document}